# Fast Stochastic Methods for Nonsmooth Nonconvex Optimization

Sashank J. Reddi
sjakkamr@cs.cmu.edu
Carnegie Mellon University

Suvrit Sra
suvrit@mit.edu
Massachusetts Institute of Technology

Barnabás Póczós
bapoczos@cs.cmu.edu
Carnegie Mellon University

Alex Smola
alex@smola.org
Carnegie Mellon University

**Abstract**

We analyze stochastic algorithms for optimizing nonconvex, nonsmooth finite-sum problems, where the nonconvex part is smooth and the nonsmooth part is convex. Surprisingly, unlike the smooth case, our knowledge of this fundamental problem is very limited. For example, it is not known whether the proximal stochastic gradient method with *constant* minibatch converges to a stationary point. To tackle this issue, we develop fast stochastic algorithms that provably converge to a stationary point for constant minibatches. Furthermore, using a variant of these algorithms, we show provably faster convergence than batch proximal gradient descent. Finally, we prove global linear convergence rate for an interesting subclass of nonsmooth nonconvex functions, that subsumes several recent works. This paper builds upon our recent series of papers on fast stochastic methods for smooth nonconvex optimization [22, 23], with a novel analysis for nonconvex and nonsmooth functions.

## 1 Introduction

We study nonconvex, nonsmooth, finite-sum optimization problems of the form

$$\min_{x \in \mathbb{R}^d} \quad F(x) := f(x) + h(x), \text{ where } f(x) := \frac{1}{n}\sum_{i=1}^n f_i(x), \tag{1}$$

where each $f_i : \mathbb{R}^d \to \mathbb{R}$ is smooth (possibly nonconvex) for all $i \in \{1, \ldots, n\} \triangleq [n]$, while $h : \mathbb{R}^d \to \mathbb{R}$ is nonsmooth but convex and relatively simple.

Such finite-sum optimization problems are fundamental to machine learning, where they typically arise within the spectrum of regularized empirical risk minimization. While there has been extensive research in solving nonsmooth *convex* finite-sum problems (i.e., each $f_i$ is convex for $i \in [n]$) [16, 4, 32], our understanding of their nonsmooth *nonconvex* counterpart is surprisingly limited—despite the widespread use and importance of nonconvex models. We focus, therefore, on fast stochastic methods for solving nonconvex, nonsmooth, finite-sum problems.

A popular approach to handle nonsmoothness is via proximal operators [14, 25]. For a proper closed convex function $h$, the *proximal operator* is defined as

$$\text{prox}_{\eta h}(x) := \operatorname*{argmin}_{y \in \mathbb{R}^d}\left(h(y) + \frac{1}{2\eta}\|y - x\|^2\right), \qquad \text{for } \eta > 0. \tag{2}$$

The power of proximal operators lies in how they generalize projections—indeed, if $h$ is the *indicator function* $\mathcal{I}_C(x)$ of a closed convex set $C$, then $\text{prox}_{\mathcal{I}_C}(x) \equiv \text{proj}_C(x) \equiv \text{argmin}_{y \in C} \|y - x\|$.

Throughout this paper, we assume that the proximal operator of $h$ is relatively easy to compute. This is true for many applications in machine learning and statistics including $\ell_1$ regularization, box-constraints, simplex constraints, among others [18, 2]. Specifically, we assume access to a *proximal oracle* (PO) that takes a point $x \in \mathbb{R}^d$ and returns the output of (2). To describe our complexity results more precisely we use the incremental first-order oracle (IFO).[1] For a function $f = \frac{1}{n} \sum_i f_i$, an IFO takes an index $i \in [n]$ and a point $x \in \mathbb{R}^d$, and returns the pair $(f_i(x), \nabla f_i(x))$.

A standard (batch) method for solving (1) is the proximal-gradient method (PROXGD) [13], first studied for nonconvex problems in [5]. This method performs the following iteration:

$$x^{t+1} = \text{prox}_{\eta h}(x^t - \eta \nabla f(x^t)), \qquad t = 0, 1, \ldots, \tag{3}$$

where $\eta > 0$ is the step size. The following non-asymptotic rate of convergence result for the proximal gradient method was proved recently.

**Theorem (Informal).** [7]: *The number of IFO and PO calls made by the proximal gradient method (3) to reach $\epsilon$ close to a stationary point is $O(n/\epsilon)$ and $O(1/\epsilon)$ respectively.*

We refer the readers to [7] for more details. The key point to note here is that the IFO complexity of (3) is $O(n/\epsilon)$. This is due to the fact that a full gradient $\nabla f$ needs to computed at each iteration of (3), thus, entailing $n$ IFO calls at each iteration. When $n$ is large, this per iteration cost is very expensive, and hence often results in slow convergence. A more practical approach is offered by the proximal stochastic gradient (PROXSGD) method, which performs the iteration

$$x^{t+1} = \text{prox}_{\eta_t h} \left( x^t - \frac{\eta_t}{|I_t|} \sum_{i \in I_t} \nabla f_i(x^t) \right), \quad t = 0, 1, \ldots, \tag{4}$$

where $I_t$ (referred to as minibatch) is a randomly chosen set (with replacement) from $[n]$ and $\eta_t$ is a step size. Non-asymptotic convergence of PROXSGD was also shown recently, as noted below.

**Theorem (Informal).** [7]: *The number of IFO and PO calls made by PROXSGD, i.e., iteration (4), to reach $\epsilon$ close to a stationary point is $O(1/\epsilon^2)$ and $O(1/\epsilon)$ respectively. For achieving this convergence, we need batch sizes $|I_t|$ that increase and step sizes $\eta_t$ that decrease with $1/\epsilon$.*

Notice that the PO complexity of PROXSGD is similar to proximal gradient, but its IFO complexity is independent of $n$; though this benefit comes at the cost of an extra $1/\epsilon$ factor. Furthermore, the step size must decrease with $1/\epsilon$ (or alternatively decay with the number of iterations of the algorithm). The same two aspects are also seen for *convex* stochastic gradient, in both the smooth and proximal versions. However, in the nonconvex setting there is a key third and more important aspect: *the minibatch size $|I_t|$ increases with $1/\epsilon$*.

To understand this aspect, consider the case of $|I_t|$ being a constant (independent of both $n$ and $\epsilon$), typically the choice used in practice. In this case, the above PROXSGD convergence result no longer holds and it is *not* clear if PROXSGD even converges to a stationary point at all. To clarify, a decreasing step size $\eta_t$ trivially ensures convergence as $t \to \infty$, but the limiting point is not necessarily stationary. On the other hand, increasing $|I_t|$ with $1/\epsilon$ can easily lead to $|I_t| \geq n$ for reasonably small $\epsilon$, which effectively reduces the algorithm to (batch) proximal gradient.

This dismal news does not apply to the convex setting, where convergence (in expectation) to an optimal solution has been shown for PROXSGD and its variants using constant minibatch sizes $|I_t|$ [24, 3]. Furthermore, this problem does not afflict the smooth nonconvex case ($h \equiv 0$), where convergence with constant minibatches is ensured [6, 22, 23]. Hence, there appears to be a fundamental gap in our understanding of stochastic methods for *nonsmooth nonconvex* problems.

---

[1] Introduced in [1] to study lower bounds of deterministic algorithms for convex finite-sum problems.

Given the ubiquity of nonconvex models in machine learning and statistics, it is important to bridge this gap. To this end, we study fast stochastic methods for tackling nonsmooth nonconvex problems with guaranteed convergence for constant minibatches, and faster convergence with minibatches independent of $1/\epsilon$.

**Main Contributions**

We state our main contributions below and list the key complexity results in Table 1.

- We analyze nonconvex proximal versions of the recently proposed stochastic algorithms SVRG and SAGA [8, 4, 32], hereafter referred to as PROXSVRG and PROXSAGA, respectively. We show convergence of these algorithms with constant minibatches. To the best of our knowledge, this is the first work to present non-asymptotic convergence rates for stochastic methods that apply to *nonsmooth nonconvex* problems with *constant* (hence more realistic) minibatches.

- We show that by carefully choosing the minibatch size (to be sublinearly dependent on $n$ but still independent of $1/\epsilon$), we can achieve provably faster convergence than both proximal gradient and proximal stochastic gradient. We are not aware of any earlier results on stochastic methods for the general *nonsmooth nonconvex* problem that have faster convergence than proximal gradient.

- We study a nonconvex subclass of (1) based on the proximal extension of Polyak-Łojasiewicz inequality [9]. We show linear convergence of PROXSVRG and PROXSAGA to the optimal solution for this subclass. This includes the recent results proved in [27, 34] as special cases. Ours is the first *stochastic* method with provable global linear convergence for this subclass of problems.

## 1.1 Related Work

The literature on finite-sum problems is vast; so we summarize only a few closely related works. Convex instances of (1) have been long studied [19, 15, 3] and are fairly well-understood. Remarkable recent progress for smooth convex instances of (1) is the creation of variance reduced (VR) stochastic methods [26, 8, 4, 28]. Nonsmooth proximal VR stochastic algorithms are studied in [32, 4] where faster convergence rates for both strongly convex and non-strongly convex cases are proved. Asynchronous VR frameworks are developed in [21]; lower-bounds are studied in [1, 10].

In contrast, nonconvex instances of (1) are much less understood. Stochastic gradient for smooth nonconvex problems is analyzed in [6], and only very recently, convergence results for VR stochastic methods for smooth nonconvex problems were obtained in [22, 23, 33]. In [11], the authors consider a VR nonconvex setting different from ours, namely, where the loss is (essentially strongly) convex but hard thresholding is used. We build upon [22, 23], and focus on handling nonsmooth convex regularizers ($h \not\equiv 0$ in (1)). Incremental proximal gradient methods for this class were also considered in [31] but only asymptotic convergence was shown. The first analysis of a projection version of nonconvex SVRG is due to [29], who considers the special problem of PCA; see also the follow-up work [30]. Perhaps, the closest to our work is [7], where convergence of minibatch nonconvex PROXSGD method is studied. However, typical to the stochastic gradient method, the convergence is slow; moreover, no convergence for constant minibatches is provided.

## 2 Preliminaries

We assume that the function $h(x)$ in (1) is lower semi-continuous (lsc) and convex. Furthermore, we also assume that its domain $\text{dom}(h) = \{x \in \mathbb{R}^d | h(x) < +\infty\}$ is closed. We say $f$ is *L-smooth* if there is a constant $L$ such that

$$\|\nabla f(x) - \nabla f(y)\| \leq L\|x - y\|, \quad \forall\ x, y \in \mathbb{R}^d.$$

| Algorithm | IFO | PO | IFO (PL) | PO (PL) | Constant minibatch? |
|---|---|---|---|---|---|
| ProxSgd | $O\left(1/\epsilon^2\right)$ | $O\left(1/\epsilon\right)$ | $O\left(1/\epsilon^2\right)$ | $O\left(1/\epsilon\right)$ | ? |
| ProxGD | $O\left(n/\epsilon\right)$ | $O\left(1/\epsilon\right)$ | $O\left(n\kappa\log(1/\epsilon)\right)$ | $O\left(\kappa\log(1/\epsilon)\right)$ | — |
| ProxSvrg | $O(n + (n^{2/3}/\epsilon))$ | $O(1/\epsilon)$ | $O((n + \kappa n^{2/3})\log(1/\epsilon))$ | $O(\kappa\log(1/\epsilon))$ | ✓ |
| ProxSaga | $O(n + (n^{2/3}/\epsilon))$ | $O(1/\epsilon)$ | $O((n + \kappa n^{2/3})\log(1/\epsilon))$ | $O(\kappa\log(1/\epsilon))$ | ✓ |

Table 1: Table comparing the *best* IFO and PO complexity of different algorithms discussed in the paper. The complexity is measured in terms of the number of oracle calls required to achieve an $\epsilon$-accurate solution. The IFO (PL) and PO (PL) represents the IFO and PO complexity of PL functions (see Section 5 for a formal definiton). The results marked in red are the contributions of this paper. In the table, "constant minibatch" indicates whether stochastic algorithm converges using a constant minibatch size. To the best of our knowledge, it is not known if ProxSgd converges on using constant minibatches for nonconvex nonsmooth optimization. Also, we are not aware of any specific convergence results for ProxSgd in the context of PL functions.

Throughout, we assume that the functions $f_i$ in (1) are $L$-smooth, so that $\|\nabla f_i(x) - \nabla f_i(y)\| \leq L\|x - y\|$ for all $i \in [n]$. Such an assumption is typical in the analysis of first-order methods.

One crucial aspect of the analysis for nonsmooth nonconvex problems is the convergence criterion. For convex problems, typically the optimality gap $F(x) - F(x^*)$ is used as a criterion. It is unreasonable to use such a criterion for general nonconvex problems due to their intractability. For smooth nonconvex problems (i.e., $h \equiv 0$), it is typical to measure stationarity, e.g., using $\|\nabla F\|$. This cannot be used for nonsmooth problems, but a fitting alternative is the *gradient mapping*[2] [17]:

$$\mathcal{G}_\eta(x) := \tfrac{1}{\eta}[x - \text{prox}_{\eta h}(x - \eta\nabla f(x))]. \tag{5}$$

When $h \equiv 0$ this mapping reduces to $\mathcal{G}_\eta(x) = \nabla f(x) = \nabla F(x)$, the gradient of function $F$ at $x$. We analyze our algorithms using the gradient mapping (5) as described more precisely below.

**Definition 1.** *A point $x$ output by stochastic iterative algorithm for solving (1) is called an $\epsilon$-accurate solution, if $\mathbb{E}[\|\mathcal{G}_\eta(x)\|^2] \leq \epsilon$ for some $\eta > 0$.*

Our goal is to obtain *efficient* algorithms for achieving an $\epsilon$-accurate solution, where efficiency is measured using IFO and PO complexity as functions of $1/\epsilon$ and $n$.

## 3 Algorithms

We focus on two algorithms: (a) proximal Svrg (ProxSvrg) and (b) proximal Saga (ProxSaga).

### 3.1 Nonconvex Proximal SVRG

We first consider a variant of ProxSvrg [32]; pseudocode of this variant is stated in Algorithm 1. When $F$ is strongly convex, Svrg attains linear convergence rate as opposed to sublinear convergence of Sgd [8]. Note that, while Svrg is typically stated with $b = 1$, we use its minibatch variant with batch size $b$. The specific reasons for using such a variant will become clear during the analysis.

While some other algorithms have been proposed for reducing the variance in the stochastic gradients, Svrg is particularly attractive because of its low memory requirement; it requires just $O(d)$ extra memory in comparison to SGD for storing the average gradient ($g^s$ in Algorithm 1), while algorithms like SAG and SAGA incur $O(nd)$ storage cost. In addition to its strong theoretical results, SVRG is known to outperform SGD empirically while being more robust to selection of step size. For convex problems, PROXSVRG is known to inherit these advantages of SVRG [32].

---

[2]This mapping has also been used in the analysis of nonconvex proximal methods in [31, 6, 7].

We now present our analysis of nonconvex PROXSVRG, starting with a result for batch size $b = 1$.

**Theorem 1.** *Let $b = 1$ in Algorithm 1. Let $\eta = 1/(3Ln)$, $m = n$ and $T$ be a multiple of $m$. Then the output $x_a$ of Algorithm 1 satisfies the following bound:*

$$\mathbb{E}[\|\mathcal{G}_\eta(x_a)\|^2] \leq \frac{18Ln^2}{3n-2}\left(\frac{F(x^0) - F(x^*)}{T}\right),$$

*where $x^*$ is an optimal solution of* (1).

Theorem 1 shows that PROXSVRG converges for constant minibatches of size $b = 1$. This result is in strong contrast to PROXSGD whose convergence with constant minibatches is still unknown. However, the result delivered by Theorem 1 is *not* stronger than that of PROXGD. The following corollary to Theorem 1 highlights this point.

**Corollary 1.** *To obtain an $\epsilon$-accurate solution, with $b = 1$ and parameters from Theorem 1, the IFO and PO complexities of Algorithm 1 are $O(n/\epsilon)$ and $O(n/\epsilon)$, respectively.*

Corollary 1 follows upon noting that each inner iteration (Step 7) of Algorithm 1 has an effective IFO comlexity of $O(1)$ since $m = n$. This IFO complexity includes the IFO calls for calculating the average gradient at the end of each epoch. Furthermore, each inner iteration also invokes the proximal oracle, whereby the PO complexity is also $O(n/\epsilon)$. While the IFO complexity of constant minibatch PROXSVRG is same as PROXGD, we see that its PO complexity is much worse. This is due to the fact that $n$ IFO calls correspond to one PO call in PROXGD, while one IFO call in PROXSVRG corresponds to one PO call. Consequently, we do not gain any theoretical advantage by using constant minibatch PROXSVRG over PROXGD.

The key question is therefore: *can we modify the algorithm to obtain better theoretical guarantees?* To answer this question, we prove the following main convergence result. For the ease of theoretical exposition, we assume $n^{2/3}$ to be an integer. This is only for convenience in stating our theoretical results and all the results in the paper hold for the general case.

**Theorem 2.** *Suppose $b = n^{2/3}$ in Algorithm 1. Let $\eta = 1/(3L)$, $m = \lfloor n^{1/3} \rfloor$ and $T$ is a multiple of $m$. Then for the output $x_a$ of Algorithm 1, we have:*

$$\mathbb{E}[\|\mathcal{G}_\eta(x_a)\|^2] \leq \frac{18L(F(x^0) - F(x^*))}{T},$$

*where $x^*$ is an optimal solution to* (1).

Rewriting Theorem 2 in terms of the IFO and PO complexity, we obtain the following corollary.

**Corollary 2.** *Let $b = n^{2/3}$ and set parameters as in Theorem 2. Then, to obtain an $\epsilon$-accurate solution, the IFO and PO complexities of Algorithm 1 are $O(n + n^{2/3}/\epsilon)$ and $O(1/\epsilon)$, respectively.*

The above corollary is due to the following observations. From Theorem 2, it can be seen that the total number of inner iterations (across all epochs) of Algorithm 1 to obtain an $\epsilon$-accurate solution is $O(1/\epsilon)$. Since each inner iteration of Algorithm 2 involves a call to the PO, we obtain a PO complexity of $O(1/\epsilon)$. Further, since $b = n^{2/3}$ IFO calls are made at each inner iteration, we obtain a net IFO complexity of $O(n^{2/3}/\epsilon)$. Adding the IFO calls for the calculation of the average gradient (and noting that $T$ is a multiple of $m$), we obtain the desired result. A noteworthy aspect of Corollary 2 is that its PO complexity matches PROXGD, but its IFO complexity is significantly decreased to $O(n + n^{2/3}/\epsilon)$ as opposed to $O(n/\epsilon)$ in PROXGD.

---
**Algorithm 1:** Nonconvex ProxSvrg $(x^0, T, m, b, \eta)$
---
1: **Input:** $\tilde{x}^0 = x_m^0 = x^0 \in \mathbb{R}^d$, epoch length $m$, step sizes $\eta > 0$, $S = \lceil T/m \rceil$
2: **for** $s = 0$ **to** $S - 1$ **do**
3: $\quad x_0^{s+1} = x_m^s$
4: $\quad g^{s+1} = \frac{1}{n} \sum_{i=1}^n \nabla f_i(\tilde{x}^s)$
5: $\quad$ **for** $t = 0$ **to** $m - 1$ **do**
6: $\quad\quad$ Uniformly randomly pick $I_t \subset \{1, \ldots, n\}$ (with replacement) such that $|I_t| = b$
7: $\quad\quad v_t^{s+1} = \frac{1}{b} \sum_{i_t \in I_t} (\nabla f_{i_t}(x_t^{s+1}) - \nabla f_{i_t}(\tilde{x}^s)) + g^{s+1}$
8: $\quad\quad x_{t+1}^{s+1} = \text{prox}_{\eta h}(x_t^{s+1} - \eta v_t^{s+1})$
9: $\quad$ **end for**
10: $\quad \tilde{x}^{s+1} = x_m^{s+1}$
11: **end for**
12: **Output:** Iterate $x_a$ chosen uniformly at random from $\{\{x_t^{s+1}\}_{t=0}^{m-1}\}_{s=0}^{S-1}$.
---

## 3.2 Nonconvex Proximal SAGA

In the previous section, we investigated ProxSvrg for solving (1). Note that ProxSvrg is not a fully "incremental" algorithm since it requires calculation of the full gradient once per epoch. An alternative to ProxSvrg is the algorithm proposed in [4] (popularly referred to as SAGA). We build upon the work of [4] to develop ProxSaga, a nonconvex proximal variant of Saga.

The pseudocode for ProxSaga is presented in Algorithm 2. The key difference between Algorithm 1 and 2 is that ProxSaga, unlike ProxSvrg, avoids computation of the full gradient. Instead, it maintains an average gradient vector $g^t$, which changes at each iteration (refer to [21]). However, such a strategy entails additional storage costs. In particular, for implementing Algorithm 2, we must store the gradients $\{\nabla f_i(\alpha_i^t)\}_{i=1}^n$, which in general can cost $O(nd)$ in storage. Nevertheless, in some scenarios common to machine learning (see [4]), one can reduce the storage requirements to $O(n)$. Whenever such an implementation of ProxSaga is possible, it can perform similar to or even better than ProxSvrg [4]; hence, in addition to theoretical interest, it is of significant practical value.

We remark that ProxSaga in Algorithm 2 differs slightly from [4]. In particular, it uses minibatches where two sets $I_t, J_t$ are sampled at each iteration as opposed to one in [4]. This is mainly for the ease of theoretical analysis.

We prove that as in the convex case, nonconvex ProxSvrg and ProxSaga share similar theoretical guarantees. In particular, our first result for ProxSaga is a counterpart to Theorem 1 for ProxSvrg.

**Theorem 3.** *Suppose $b = 1$ in Algorithm 2. Let $\eta = 1/(5Ln)$. Then for the output $x_a$ of Algorithm 2 after $T$ iterations, the following stationarity bound holds:*

$$\mathbb{E}[\|\mathcal{G}_\eta(x_a)\|^2] \leq \frac{50Ln^2}{5n-2} \frac{F(x^0) - F(x^*)}{T},$$

*where $x^*$ is an optimal solution of (1).*

Theorem 3 immediately leads to the following corollary.

**Corollary 3.** *The IFO and PO complexity of Algorithm 3 for $b = 1$ and parameters specified in Theorem 3 to obtain an $\epsilon$-accurate solution are $O(n/\epsilon)$ and $O(n/\epsilon)$ respectively.*

Similar to Theorem 2 for ProxSvrg, we obtain the following main result for ProxSaga.

---

**Algorithm 2:** Nonconvex PROXSAGA $(x^0, T, b, \eta)$

1: **Input:** $x^0 \in \mathbb{R}^d$, $\alpha_i^0 = x^0$ for $i \in [n]$, step size $\eta > 0$
2: $g^0 = \frac{1}{n} \sum_{i=1}^n \nabla f_i(\alpha_i^0)$
3: **for** $t = 0$ to $T - 1$ **do**
4:     Uniformly randomly pick sets $I_t, J_t$ from $[n]$ (with replacement) such that $|I_t| = |J_t| = b$
5:     $v^t = \frac{1}{b} \sum_{i_t \in I_t} (\nabla f_{i_t}(x^t) - \nabla f_{i_t}(\alpha_{i_t}^t)) + g^t$
6:     $x^{t+1} = \text{prox}_{\eta h}(x^t - \eta v^t)$
7:     $\alpha_j^{t+1} = x^t$ for $j \in J_t$ and $\alpha_j^{t+1} = \alpha_j^t$ for $j \notin J_t$
8:     $g^{t+1} = g^t - \frac{1}{n} \sum_{j_t \in J_t} (\nabla f_{j_t}(\alpha_{j_t}^t) - \nabla f_{j_t}(\alpha_{j_t}^{t+1}))$
9: **end for**
10: **Output:** Iterate $x_a$ chosen uniformly random from $\{x^t\}_{t=0}^{T-1}$.

---

**Theorem 4.** *Suppose $b = n^{2/3}$ in Algorithm 2. Let $\eta = 1/(5L)$. Then for the output $x_a$ of Algorithm 2 after $T$ iterations, the following holds:*

$$\mathbb{E}[\|\mathcal{G}_\eta(x_a)\|^2] \leq \frac{50L(F(x^0) - F(x^*))}{3T},$$

*where $x^*$ is an optimal solution of Problem* (1).

Rewriting this result in terms of IFO and PO access, we obtain the following important corollary.

**Corollary 4.** *Let $b = n^{2/3}$ and set parameters as in Theorem 4. Then, to obtain an $\epsilon$-accurate solution, the IFO and PO complexities of Algorithm 2 are $O(n + n^{2/3}/\epsilon)$ and $O(1/\epsilon)$, respectively.*

The above result is due to Theorem 4 and because each iteration of PROXSAGA requires $O(n^{2/3})$ IFO calls. The number of PO calls is only $O(1/\epsilon)$, since make one PO call for every $n^{2/3}$ IFO calls.

**Discussion**: It is important to note the role of minibatches in Corollaries 2 and 4. Minibatches are typically used for reducing variance and promoting parallelism in stochastic methods. But unlike previous works, we use minibatches as a theoretical tool to improve convergence rates of both nonconvex PROXSVRG and PROXSAGA. In particular, by carefully selecting the minibatch size, we can improve the IFO complexity of the algorithms described in the paper from $O(n/\epsilon)$ (similar to PROXGD) to $O(n^{2/3}/\epsilon)$ (matching our results for the smooth nonconvex case in [22, 23]). Furthermore, the PO complexity is also improved in a similar manner by using the minibatch size mentioned in Theorems 2 and 4.

## 4 General Convergence Analysis

In this previous sections, for the sake of clarity, we stated convergence rates of PROXSVRG and PROXSAGA for a specific set of parameters. However, a more general analysis can be derived for these algorithms. The rationale behind the choice of parameters is Section 3 will also become clear later in this section. We have the following general convergence results for PROXSVRG and PROXSAGA.

**Theorem 5.** *Suppose $b \leq n$ in Algorithm 1. Let $T$ be a multiple of $m$ and $\eta = \rho/L$ where $\rho < 1/2$ and satisfies the following:*

$$\frac{4\rho^2 m^2}{b} + \rho \leq 1.$$

*Then for the output $x_a$ of Algorithm 1, we have:*

$$\mathbb{E}[\|\mathcal{G}_\eta(x_a)\|^2] \leq \frac{2L(F(x^0) - F(x^*))}{\rho(1 - 2\rho)T},$$

*where $x^*$ is an optimal solution to* (1).

The following result is an immediate consequence of the above result.

**Corollary 5.** *Let $b \leq n$, $\rho = 1/4$ and $m = \lfloor b^{1/2} \rfloor$ in Theorem 5. Then, to obtain an $\epsilon$-accurate solution, the IFO and PO complexities of Algorithm 1 are $O(n + n/(b^{1/2}\epsilon) + b/\epsilon)$ and $O(1/\epsilon)$, respectively.*

We observe that under the parameter setting in Corollary 5, the PO complexity is $O(1/\epsilon)$, which matches that of ProxGD. Thus, this setting is optimized for reducing the PO complexity. Furthermore, for constant minibatch $b = 1$, Corollary 5 shows that the IFO and PO complexity of ProxSvrg is $O(n/\epsilon)$ and $O(1/\epsilon)$ respectively, which is stronger than Theorem 1. However, in the setting of Corollary 5 with minibatch size $b = 1$, ProxSvrg effectively reduces to ProxGD since $m = \lfloor b^{1/2} \rfloor = 1$ and hence, not very interesting. When $b = 1$, Theorem 1 provides the convergence result of ProxSvrg for the setting that is generally used in practice where $m = n$.

For ProxSaga, we have the following general convergence results.

**Theorem 6.** *Suppose $b \leq n$ in Algorithm 2. Let $\eta = \rho/L$ where $\rho < 1/2$ and $\rho$ satisfies the following condition:*
$$\frac{16n^2\rho^2}{b^3} + \rho \leq 1.$$
*Then for the output $x_a$ of Algorithm 2, we have:*
$$\mathbb{E}[\|\mathcal{G}_\eta(x_a)\|^2] \leq \frac{2L(F(x^0) - F(x^*))}{\rho(1 - 2\rho)T},$$
*where $x^*$ is an optimal solution to (1).*

**Corollary 6.** *Let $b \leq n$, $\rho = \min\{1/5, b^{3/2}/5n\}$ in Theorem 6. Then, to obtain an $\epsilon$-accurate solution, the IFO and PO complexities of Algorithm 2 are $O(n + n/(b^{1/2}\epsilon) + b/\epsilon)$ and $O(\max\{1, n/b^{3/2}\}/\epsilon)$ respectively.*

Note that the IFO complexity of ProxSvrg (in Corollary 5) and ProxSaga (in Corollary 6) are similar; however, their PO complexities are different. It is not hard to see from Corollary 5 and 6, that the best IFO and PO complexity of both ProxSvrg and ProxSaga obtainable through these upper bounds are $O(n + n^{2/3}/\epsilon)$ and $O(1/\epsilon)$ respectively; which are precisely our main results in Section 3.

## 5 Extensions

We discuss some extensions of our approach in this section. Our first extension is to provide convergence analysis for a subclass of nonconvex functions that satisfy a specific growth condition popularly known as the Polyak-Łojasiewicz (PL) inequality. In the context of gradient descent, this inequality was proposed by Polyak in 1963 [20], who showed *global* linear convergence of gradient descent for functions that satisfy the PL inequality. Recently, in [9] the PL inequality was generalized to nonsmooth functions and used for proving linear convergence of proximal gradient. The generalization presented in [9] considers functions $F(x) = f(x) + h(x)$ that satisfy the following:

$$\mu(F(x) - F(x^*)) \leq \frac{1}{2} D_h(x, \mu), \text{ where } \mu > 0$$
$$\text{and } D_h(x, \mu) := -2\mu \min_y \left[\langle \nabla f(x), y - x \rangle + \frac{\mu}{2}\|y - x\|^2 + h(y) - h(x)\right]. \tag{6}$$

An $F$ that satisfies (6) is called a $\mu$-PL function.

When $h \equiv 0$, condition (6) reduces to the usual PL inequality. The class of $\mu$-PL functions includes several other classes as special cases. It subsumes strongly convex functions, covers $f_i(x) =$ $g(a_i^\top x)$ with only $g$ being strongly convex, and includes functions that satisfy a optimal strong convexity property [12]. Note that the $\mu$-PL functions also subsume the recently studied special case where $f_i$'s are nonconvex but their sum $f$ is strongly convex. Hence, it encapsulates the problems of [27, 34].

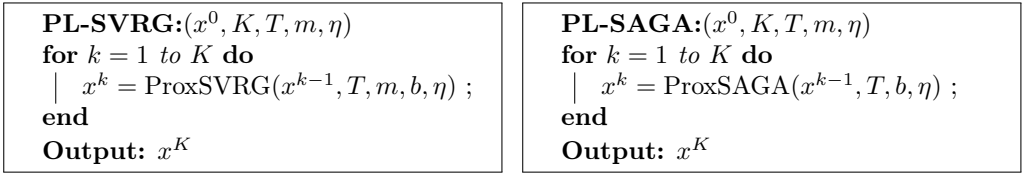

Figure 1: PROXSVRG and PROXSAGA variants for PL functions.

The algorithms in Figure 1 provide variants of PROXSVRG and PROXSAGA adapted to optimize $\mu$-PL functions. We show the following *global* linear convergence result of PL-SVRG and PL-SAGA in Figure 1 for PL functions. For simplicity, we assume $\kappa = (L/\mu) > n^{1/3}$. When $f$ is strongly convex, $\kappa$ is referred to as the condition number, in which case $\kappa > n^{1/3}$ corresponds to the high condition number regime.

**Theorem 7.** *Suppose $F$ is a $\mu$-PL function. Let $b = n^{2/3}$, $\eta = 1/5L$, $m = \lfloor n^{1/3} \rfloor$ and $T = \lceil 30\kappa \rceil$. Then for the output $x^K$ of* PL-SVRG *and* PL-SAGA *(in Figure 1), the following holds:*

$$\mathbb{E}[F(x^K) - F(x^*)] \leq \frac{[F(x^0) - F(x^*)]}{2^K},$$

*where $x^*$ is an optimal solution of* (1).

The following corollary on IFO and PO complexity of PL-SVRG and PL-SAGA is immediate.

**Corollary 7.** *When $F$ is a $\mu$-PL function, then the IFO and PO complexities of* PL-SVRG *and* PL-SAGA *with the parameters specified in Theorem 7 to obtain an $\epsilon$-accurate solution are $O((n + \kappa n^{2/3})\log(1/\epsilon))$ and $O(\kappa \log(1/\epsilon))$, respectively.*

Note that proximal gradient also has global linear convergence for PL functions, as recently shown in [9]. However, its IFO complexity is $O(\kappa n \log(1/\epsilon))$, which is much worser than that of PL-SVRG and PL-SAGA (Corollary 7).

**Other extensions:** Our results can be easily generalized to the case where non-uniform sampling is used in Algorithm 1 and Algorithm 2. This is useful when the functions $f_i$ have different Lipschitz constants.

## 6 Experiments

We present our empirical results in this section. For our experiments, we study the problem of non-negative principal component analysis (NN-PCA). More specifically, for a given set of samples $\{z_i\}_{i=1}^n$, we solve the following optimization problem:

$$\min_{\|x\|\leq 1,\ x\geq 0} -\frac{1}{2} x^\top \left( \sum_{i=1}^n z_i z_i^\top \right) x. \tag{7}$$

The problem of NN-PCA is, in general, NP-hard. This variant of the standard PCA problem can be written in the form (1) with $f_i(x) = -(x^\top z_i)^2$ for all $i \in [n]$ and $h(x) = \mathcal{I}_C(x)$ where $C$ is the convex set $\{x \in \mathbb{R}^d | \|x\| \leq 1, x \geq 0\}$. In our experiments, we compare PROXSGD with nonconvex

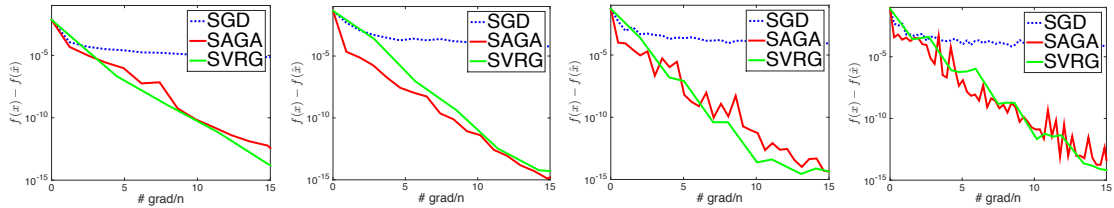

Figure 2: Non-negative principal component analysis. Performance of PROXSGD, PROXSVRG and PROXSAGA on 'rcv1' (left), 'a9a'(left-center), 'mnist' (right-center) and 'aloi' (right) datasets. Here, the y-axis is the function suboptimality i.e., $f(x) - f(\hat{x})$ where $\hat{x}$ represents the best solution obtained by running gradient descent for long time and with multiple restarts.

PROXSVRG and PROXSAGA. The choice of step size is important to PROXSGD. The step size of PROXSGD is set using the popular $t$-inverse step size choice of $\eta_t = \eta_0(1 + \eta'\lfloor t/n \rfloor)^{-1}$ where $\eta_0, \eta' > 0$. For PROXSVRG and PROXSAGA, motivated by the theoretical analysis, we use a fixed step size. The parameters of the step size in each of these methods are chosen so that the method gives the best performance on the objective value. In our experiments, we include the value $\eta' = 0$, which corresponds to PROXSGD with fixed step size. For PROXSVRG, we use the epoch length $m = n$.

We use standard machine learning datasets in LIBSVM for all our experiments [3]. The samples from each of these datasets are normalized i.e. $\|z_i\| = 1$ for all $i \in [n]$. Each of these methods is initialized by running PROXSGD for $n$ iterations. Such an initialization serves two purposes: (a) it provides a reasonably good initial point, typically beneficial for variance reduction techniques [26, 4]. (b) it provides a heuristic for calculating the initial average gradient $g^0$ [26]. In our experiments, we use minibatch size $b = 1$ in order to demonstrate the performance of the algorithms with constant minibatches.

We report the objective function value for the datasets. In particular, we report the suboptimality in objective function i.e., $f(x_t^{s+1}) - f(\hat{x})$ (for PROXSVRG) and $f(x^t) - f(\hat{x})$ (for PROXSAGA). Here $\hat{x}$ refers to the solution obtained by running proximal gradient descent for a large number of iterations and multiple random initializations. For all the algorithms, we compare the aforementioned criteria against for the number of *effective* passes through the dataset i.e., IFO complexity divided by $n$. For PROXSVRG, this includes the cost of calculating the full gradient at the end of each epoch.

Figure 2 shows the performance of PROXSGD , PROXSVRG and PROXSVRG on NN-PCA problem (see Section D of the Appendix for more experiments). It can be seen that the objective value for PROXSVRG and PROXSAGA is much lower compared to PROXSGD, suggesting faster convergence for these algorithms. We observed a significant gain consistently across all the datasets. Moreover, the selection of step size was much simpler for PROXSVRG and PROXSAGA than that for PROXSGD. We did not observe any significant difference in the performance of PROXSVRG and PROXSAGA for this particular task.

## 7 Final Discussion

In this paper, we presented fast stochastic methods for nonsmooth nonconvex optimization. In particular, by employing variance reduction techniques, we show that one can design methods that can provably perform better than PROXSGD and proximal gradient descent. Furthermore, in contrast to PROXSGD, the resulting approaches have provable convergence to a stationary point with constant minibatches; thus, bridging a fundamental gap in our knowledge of nonsmooth nonconvex problems.

We proved that with a careful selection of minibatch size, it is possible to theoretically show superior performance to proximal gradient descent. Our empirical results provide evidence for a similar conclusion even with constant minibatchs. Thus, we conclude with an important open problem of developing stochastic methods with provably better performance than proximal gradient descent with *constant minibatch* size.

---

[3]The datasets can be downloaded from https://www.csie.ntu.edu.tw/~cjlin/libsvmtools/datasets.

## Acknowledgment

We thank Guanghui Lan for very enlightening discussions.

## A  Convergence analysis for Proximal Nonconvex SVRG

The analysis requires some key lemmas which can be found in Appendix E.

### A.1  General Convergence Analysis: Proof of Theorem 5

*Proof.* We start by defining the full gradient iterate

$$\overline{x}_{t+1}^{s+1} = \text{prox}_{\eta h}(x_t^{s+1} - \eta \nabla f(x_t^{s+1})), \quad (8)$$

which is merely for our analysis, and is never actually computed. Applying Lemma 2 to (8) (with $y = \overline{x}_{t+1}^{s+1}$, $z = x_t^{s+1}$ and $d' = \nabla f(x_t^{s+1})$), and taking expectations we obtain the bound

$$\mathbb{E}[F(\overline{x}_{t+1}^{s+1})] \leq \mathbb{E}\left[F(x_t^{s+1}) + \left[\tfrac{L}{2} - \tfrac{1}{2\eta}\right]\|\overline{x}_{t+1}^{s+1} - x_t^{s+1}\|^2 - \tfrac{1}{2\eta}\|\overline{x}_{t+1}^{s+1} - x_t^{s+1}\|^2\right]. \quad (9)$$

Recall the iterates of Algorithm 1 are computed using the following update:

$$x_{t+1}^{s+1} = \text{prox}_{\eta h}(x_t^{s+1} - \eta v_t^{s+1})), \quad (10)$$

where $v_t^{s+1} = \frac{1}{b}\sum_{i_t \in I_t}(\nabla f_{i_t}(x_t^{s+1}) - \nabla f_{i_t}(\tilde{x}^s)) + g^{s+1}$ (see Algorithm 1). Applying Lemma 2 on update (10) (with $y = x_{t+1}^{s+1}$, $z = \overline{x}_{t+1}^{s+1}$ and $d' = v_t^{s+1}$) and taking expectations we obtain

$$\mathbb{E}[F(x_{t+1}^{s+1})] \leq \mathbb{E}\Big[F(\overline{x}_{t+1}^{s+1}) + \langle x_{t+1}^{s+1} - \overline{x}_{t+1}^{s+1}, \nabla f(x_t^{s+1}) - v_t^{s+1}\rangle$$
$$+ \Big[\tfrac{L}{2} - \tfrac{1}{2\eta}\Big]\|x_{t+1}^{s+1} - x_t^{s+1}\|^2 + \Big[\tfrac{L}{2} + \tfrac{1}{2\eta}\Big]\|\overline{x}_{t+1}^{s+1} - x_t^{s+1}\|^2 - \tfrac{1}{2\eta}\|x_{t+1}^{s+1} - \overline{x}_{t+1}^{s+1}\|^2\Big]. \quad (11)$$

Adding inequalities (9) and (11), we get

$$\mathbb{E}[F(x_{t+1}^{s+1})] \leq \mathbb{E}\Big[F(x_t^{s+1}) + \Big[L - \tfrac{1}{2\eta}\Big]\|\overline{x}_{t+1}^{s+1} - x_t^{s+1}\|^2 + \Big[\tfrac{L}{2} - \tfrac{1}{2\eta}\Big]\|x_{t+1}^{s+1} - x_t^{s+1}\|^2$$
$$- \tfrac{1}{2\eta}\|x_{t+1}^{s+1} - \overline{x}_{t+1}^{s+1}\|^2 + \underbrace{\langle x_{t+1}^{s+1} - \overline{x}_{t+1}^{s+1}, \nabla f(x_t^{s+1}) - v_t^{s+1}\rangle}_{T_1}\Big] \quad (12)$$

We can bound the term $T_1$ as follows:

$$\mathbb{E}[T_1] \leq \frac{1}{2\eta}\mathbb{E}\|x_{t+1}^{s+1} - \overline{x}_{t+1}^{s+1}\|^2 + \frac{\eta}{2}\mathbb{E}\|\nabla f(x_t^{s+1}) - v_t^{s+1}\|^2$$
$$\leq \frac{1}{2\eta}\mathbb{E}\|x_{t+1}^{s+1} - \overline{x}_{t+1}^{s+1}\|^2 + \frac{\eta L^2}{2b}\mathbb{E}\|x_t^{s+1} - \tilde{x}^s\|^2.$$

The first inequality follows from Cauchy-Schwarz and Young's inequality, while the second inequality is due to Lemma 3. Substituting the upper bound on $T_1$ in (12), we see that

$$\mathbb{E}[F(x_{t+1}^{s+1})] \leq \mathbb{E}\Big[F(x_t^{s+1}) + \Big[L - \tfrac{1}{2\eta}\Big]\|\overline{x}_{t+1}^{s+1} - x_t^{s+1}\|^2$$
$$+ \Big[\tfrac{L}{2} - \tfrac{1}{2\eta}\Big]\|x_{t+1}^{s+1} - x_t^{s+1}\|^2 + \tfrac{\eta L^2}{2b}\|x_t^{s+1} - \tilde{x}^s\|^2\Big]. \quad (13)$$

To further analyze (13), we set up a recursion for which we use the following Lyapunov function:

$$R_t^{s+1} := \mathbb{E}[F(x_t^{s+1}) + c_t\|x_t^{s+1} - \tilde{x}^s\|^2].$$

Introduce the quantities $c_m = 0$, and $c_t = c_{t+1}(1+\beta) + \frac{\eta L^2}{2b}$. Also, for rest of the analysis set $\beta = 1/m$. We can then bound $R_{t+1}^{s+1}$ as follows

$$R_{t+1}^{s+1} = \mathbb{E}[F(x_{t+1}^{s+1}) + c_{t+1}\|x_{t+1}^{s+1} - x_t^{s+1} + x_t^{s+1} - \tilde{x}^s\|^2]$$
$$= \mathbb{E}[F(x_{t+1}^{s+1}) + c_{t+1}(\|x_{t+1}^{s+1} - x_t^{s+1}\|^2 + \|x_t^{s+1} - \tilde{x}^s\|^2 + 2\langle x_{t+1}^{s+1} - x_t^{s+1}, x_t^{s+1} - \tilde{x}^s\rangle)]$$
$$\leq \mathbb{E}[F(x_{t+1}^{s+1}) + c_{t+1}(1+1/\beta)\|x_{t+1}^{s+1} - x_t^{s+1}\|^2 + c_{t+1}(1+\beta)\|x_t^{s+1} - \tilde{x}^s\|^2]$$
$$\leq \mathbb{E}\Big[F(x_t^{s+1}) + \Big[L - \tfrac{1}{2\eta}\Big]\|\overline{x}_{t+1}^{s+1} - x_t^{s+1}\|^2 + \Big[c_{t+1}\Big(1 + \tfrac{1}{\beta}\Big) + \tfrac{L}{2} - \tfrac{1}{2\eta}\Big]\|x_{t+1}^{s+1} - x_t^{s+1}\|^2$$
$$+ \Big[c_{t+1}(1+\beta) + \tfrac{\eta L^2}{2b}\Big]\|x_t^{s+1} - \tilde{x}^s\|^2\Big] \quad (14)$$
$$\leq \mathbb{E}\Big[F(x_t^{s+1}) + \Big[L - \tfrac{1}{2\eta}\Big]\|\overline{x}_{t+1}^{s+1} - x_t^{s+1}\|^2 + \Big[c_{t+1}(1+\beta) + \tfrac{\eta L^2}{2b}\Big]\|x_t^{s+1} - \tilde{x}^s\|^2\Big]$$
$$= R_t^{s+1} + \Big[L - \tfrac{1}{2\eta}\Big]\mathbb{E}\|\overline{x}_{t+1}^{s+1} - x_t^{s+1}\|^2. \quad (15)$$

The first inequality follows from Cauchy-Schwarz and Young's inequality. The second inequality is due to the bound (13), while the final equality is due to the definition of the Lyapunov function $R_t^{s+1}$. The third inequality holds because the sequence of values $c_t$ satisfies the following bound:

$$c_{t+1}\left(1 + \frac{1}{\beta}\right) + \frac{L}{2} \leq \frac{1}{2\eta}. \quad (16)$$

To verify (16), first observe that $c_m = 0$ and $c_t = c_{t+1}(1+\beta) + \frac{\eta L^2}{2b}$. Recursing on $t$, we thus obtain

$$c_t = \frac{\eta L^2}{2b} \frac{(1+\beta)^{m-t}-1}{\beta} = \frac{\rho L m}{2b}\left(\left(1+\frac{1}{m}\right)^{m-t}-1\right)$$

$$\leq \frac{\rho L m}{2b}(e-1) \leq \frac{\rho L m}{b},$$

where the first equality is due to the definition of $\eta$ and $\beta$. The first inequality follows upon noting that (i) $\lim_{l\to+\infty}(1+1/l)^l = e$ and (ii) $(1+1/l)^l$ is an increasing function for $l > 0$ (here $e$ is Euler's number). It follows that

$$c_{t+1}\left(1+\frac{1}{\beta}\right) + \frac{L}{2} \leq \frac{\rho L m}{b}(1+m) + \frac{L}{2}$$

$$\leq \frac{2\rho L m^2}{b} + \frac{L}{2} \leq \frac{L}{2\rho} = \frac{1}{2\eta},$$

where the second inequality uses $m \geq 1$. The third inequality uses the condition that

$$\frac{4\rho^2 m^2}{b} + \rho \leq 1.$$

Hence, inequality (16) follows. Now, adding (15) across all the iterations in epoch $s+1$ and then telescoping sums, we get

$$R_m^{s+1} \leq R_0^{s+1} + \sum_{t=0}^{m-1}\left[L - \frac{1}{2\eta}\right]\mathbb{E}\|\overline{x}_{t+1}^{s+1} - x_t^{s+1}\|^2. \tag{17}$$

Since $c_m = 0$ and from the definition of $\tilde{x}^{s+1}$, it follows that $R_m^{s+1} = \mathbb{E}[F(x_m^{s+1})] = \mathbb{E}[F(\tilde{x}^{s+1})]$. Furthermore, $R_0^{s+1} = \mathbb{E}[F(x_0^{s+1})] = \mathbb{E}[F(\tilde{x}^s)]$. This is due to the fact that $x_0^{s+1} = \tilde{x}^s$. Therefore, using the above reasoning in inequality (17), we have

$$\mathbb{E}[F(\tilde{x}^{s+1})] \leq \mathbb{E}[F(\tilde{x}^s)] + \sum_{t=0}^{m-1}\left[L - \frac{1}{2\eta}\right]\mathbb{E}\|\overline{x}_{t+1}^{s+1} - x_t^{s+1}\|^2. \tag{18}$$

Adding (18) across all the epochs and rearranging the terms, we obtain the bound

$$\sum_{s=0}^{S}\sum_{t=0}^{m-1}\left[\frac{1}{2\eta} - L\right]\mathbb{E}\|\overline{x}_{t+1}^{s+1} - x_t^{s+1}\|^2 \leq F(x^0) - \mathbb{E}[F(\tilde{x}^S)] \leq F(x^0) - F(x^*), \tag{19}$$

where the second inequality follows from the optimality of $x^*$.

Recall that in our notation

$$\mathcal{G}_\eta(x_t^{s+1}) = \tfrac{1}{\eta}[x_t^{s+1} - \mathrm{prox}_{\eta h}(x_t^{s+1} - \eta \nabla f(x_t^{s+1}))] = \tfrac{1}{\eta}[x_t^{s+1} - \overline{x}_{t+1}^{s+1}].$$

Using this relationship in (19) we see that

$$\sum_{s=0}^{S}\sum_{t=0}^{m-1}\left[\tfrac{1}{2\eta} - L\right]\eta^2 \mathbb{E}\|\mathcal{G}_\eta(x_t^{s+1})\|^2 \leq F(x^0) - F(x^*). \tag{20}$$

Now using the definition of $x_a$ from Algorithm 1 and simplifying we obtain the desired result. □

## Proof of Theorem 1

*Proof.* The proof follows from Theorem 5 with $b = 1$ and the parameters used in the theorem statement. □

## Proof of Theorem 2

*Proof.* The proof follows from Theorem 5 with $b = n^{2/3}$ and the parameters used in the theorem statement. □

## B  Convergence analysis for Nonconvex Proximal SAGA

### B.1  General Convergence Analysis: Proof of Theorem 6

*Proof.* We introduce the full-gradient iterate (as before, only for the analysis)

$$\overline{x}^{t+1} = \text{prox}_{\eta h}(x^t - \eta \nabla f(x^t)), \tag{21}$$

and recall that PROXSAGA iterations compute the update

$$x^{t+1} = \text{prox}_{\eta h}(x^t - \eta v^t),$$

where $v^t = \frac{1}{b} \sum_{i_t \in I_t} \left( \nabla f_{i_t}(x^t) - \nabla f_{i_t}(\alpha_{i_t}^t) \right) + g^t$. Now, using the same argument as in Theorem 5 until inequality (12), we obtain the following

$$\mathbb{E}[F(x^{t+1})] \leq \mathbb{E}\Big[F(x^t) + \left[L - \tfrac{1}{2\eta}\right] \|\overline{x}^{t+1} - x^t\|^2 + \left[\tfrac{L}{2} - \tfrac{1}{2\eta}\right] \|x^{t+1} - x^t\|^2$$
$$- \tfrac{1}{2\eta}\|x^{t+1} - \overline{x}^{t+1}\|^2 + \underbrace{\langle x^{t+1} - \overline{x}^{t+1}, \nabla f(x^t) - v^t \rangle}_{T_2}\Big]. \tag{22}$$

The term $T_2$ in (22) can be bound as follows:

$$\mathbb{E}[T_2] \leq \frac{1}{2\eta} \mathbb{E}\|x^{t+1} - \overline{x}^{t+1}\|^2 + \frac{\eta}{2} \mathbb{E}\|\nabla f(x^t) - v^t\|^2$$
$$\leq \frac{1}{2\eta} \mathbb{E}\|x^{t+1} - \overline{x}^{t+1}\|^2 + \frac{\eta L^2}{2nb} \sum_{i=1}^n \mathbb{E}\|x^t - \alpha_i^t\|^2.$$

The inequality follows from Cauchy-Schwarz and Young's inequality. The second inequality is due to Lemma 4. Substituting the upper bound on $T_2$ in inequality (22), we have

$$\mathbb{E}[F(x^{t+1})] \leq \mathbb{E}\Big[F(x^t) + \left[L - \tfrac{1}{2\eta}\right] \|\overline{x}^{t+1} - x^t\|^2$$
$$+ \left[\tfrac{L}{2} - \tfrac{1}{2\eta}\right] \|x^{t+1} - x^t\|^2 + \frac{\eta L^2}{2nb} \sum_{i=1}^n \|x^t - \alpha_i^t\|^2\Big]. \tag{23}$$

For further analysis, we require the following Lyapunov function:

$$R_t := \mathbb{E}\left[F(x^t) + \frac{c_t}{n} \sum_{i=1}^n \|x^t - \alpha_i^t\|^2\right].$$

Moreover, for the rest of the analysis we set $\beta = b/4n$. We use $p$ to denote the probability $1-(1-1/n)^b$ of an index $i$ being in $J_t$. Observe that we can bound $p$ from below as

$$p = 1 - \left(1 - \tfrac{1}{n}\right)^b \geq 1 - \tfrac{1}{1+(b/n)} = \tfrac{b/n}{1+b/n} \geq \tfrac{b}{2n}, \tag{24}$$

where the first inequality follows from $(1-y)^r \leq 1/(1+ry)$ (which holds for $y \in [0,1]$ and $r \geq 1$), while the second inequality holds because $b \leq n$.

We now obtain a recursive bound on $R_{t+1}$ as follows

$$R_{t+1} = \mathbb{E}[F(x^{t+1}) + \frac{c_{t+1}}{n}\sum_{i=1}^{n}\|x^{t+1} - \alpha_i^{t+1}\|^2]$$

$$= \mathbb{E}\big[F(x^{t+1}) + \frac{c_{t+1}p}{n}\sum_{i=1}^{n}\|x^{t+1} - x^t\|^2 + \frac{c_{t+1}(1-p)}{n}\sum_{i=1}^{n}\|x^{t+1} - \alpha_i^t\|^2\big]$$

$$= \mathbb{E}\big[F(x^{t+1}) + c_{t+1}p\|x^{t+1} - x^t\|^2$$
$$+ \frac{c_{t+1}(1-p)}{n}\sum_{i=1}^{n}(\|x^{t+1} - x^t\|^2 + \|x^t - \alpha_i^t\|^2 + 2\langle x^{t+1} - x^t, x^t - \alpha_i^t\rangle)\big]$$

$$\leq \mathbb{E}\big[F(x^{t+1}) + c_{t+1}\big(1 + \frac{1-p}{\beta}\big)\|x^{t+1} - x^t\|^2 + \frac{c_{t+1}(1+\beta)(1-p)}{n}\sum_{i=1}^{n}\|x^t - \alpha_i^t\|^2\big]$$

$$\leq \mathbb{E}\bigg[F(x^t) + \big[L - \tfrac{1}{2\eta}\big]\|\overline{x}^{t+1} - x^t\|^2 + \big[c_{t+1}\big(1 + \tfrac{1-p}{\beta}\big) + \tfrac{L}{2} - \tfrac{1}{2\eta}\big]\|x^{t+1} - x^t\|^2$$
$$+ \big[\tfrac{c_{t+1}(1+\beta)(1-p)}{n} + \tfrac{\eta L^2}{2nb}\big]\sum_{i=1}^{n}\|x^t - \alpha_i^t\|^2\bigg] \tag{25}$$

$$\leq \mathbb{E}\bigg[F(x^t) + \big[L - \tfrac{1}{2\eta}\big]\|\overline{x}^{t+1} - x^t\|^2 + \big[\tfrac{c_{t+1}(1+\beta)(1-p)}{n} + \tfrac{\eta L^2}{2nb}\big]\sum_{i=1}^{n}\|x_t - \alpha_i^t\|^2\bigg]$$

$$= R_t + \big[L - \tfrac{1}{2\eta}\big]\mathbb{E}\|\overline{x}^{t+1} - x^t\|^2. \tag{26}$$

The equality in the second line follows how $\alpha_i^{t+1}$ is chosen in Algorithm 2. In particular, from noting that each index in $J_t$ is drawn uniformly randomly and independently from $[n]$. The first inequality follows from Cauchy-Schwarz and Young's inequality. The second inequality uses the bound (23). The final equality is due to the definition of the Lyapunov function $R_t$, wherein we also use

$$c_t = \big[c_{t+1}(1+\beta)(1-p) + \frac{\eta L^2}{2b}\big]. \tag{27}$$

The third inequality requires a brief explanation. It follows upon observing that

$$c_{t+1}\left(1 + \frac{1-p}{\beta}\right) + \frac{L}{2} \leq \frac{1}{2\eta}. \tag{28}$$

To see why (28) holds, first observe that $c_T = 0$, and then use (27) to show that

$$c_t \leq c_{t+1}(1-\theta) + \frac{\eta L^2}{2b},$$

where $\theta = (b/2n) - \beta = b/4n$. The above inequality is elementary, since $(1+\beta)(1-p) \leq 1 - p + \beta \leq (1-\theta)$ and because $p \geq (b/2n)$ as noted in (24). Recursing on $t$, we thus obtain

$$c_t \leq \frac{\eta L^2}{2b}\frac{1-(1-\theta)^{T-t}}{\theta} \leq \frac{2n\rho L}{b^2}, \tag{29}$$

for all $t \in \{0, \ldots, T-1\}$, which holds due to the definition of $\eta$ and $\theta$. We now use inequality (29)

to bound the left hand side of (28) as follows

$$c_{t+1}\left(1+\frac{1-p}{\beta}\right)+\frac{L}{2} \leq \frac{2n\rho L}{b^2}\left(1+\frac{2(2n-b)}{b}\right)+\frac{L}{2}$$
$$=\frac{2n\rho L}{b^2}\left[\frac{4n}{b}-1\right]+\frac{L}{2}$$
$$\leq \frac{L}{2\rho}=\frac{1}{2\eta}.$$

The first inequality uses the bound (24), while the third inequality uses the following condition on $\rho$:

$$\frac{16n^2\rho^2}{b^3}+\rho \leq 1.$$

Thus, inequality (28) holds. Adding the bound (26) across all the iterations and then using telescoping sums, we get

$$R_T \leq R_0 + \sum_{t=0}^{T-1}\left[L-\frac{1}{2\eta}\right]\mathbb{E}\|\overline{x}^{t+1}-x^t\|^2. \qquad (30)$$

Since $c_T = 0$, we observe that $R_T = \mathbb{E}[F(x^T)]$. Furthermore, since $\alpha_i^0 = x^0$ for all $i \in [n]$, we conclude that $R_0 = \mathbb{E}[F(x^0)]$. Therefore, we can rewrite (30) to obtain

$$\mathbb{E}[F(x^T)] \leq F(x^0) + \sum_{t=0}^{T-1}\left[L-\frac{1}{2\eta}\right]\mathbb{E}\|\overline{x}^{t+1}-x^t\|^2.$$

Rearranging, and using optimality of $x^*$, this leads to the bound

$$\sum_{t=0}^{T-1}\left[\frac{1}{2\eta}-L\right]\mathbb{E}\|\overline{x}^{t+1}-x^t\|^2 \leq F(x^0)-\mathbb{E}[F(x^T)] \leq F(x^0)-F(x^*).$$

Now recall the relationship

$$\mathcal{G}_\eta(x^t) = \tfrac{1}{\eta}[x^t - \text{prox}_{\eta h}(x^t - \eta \nabla f(x^t))] = \tfrac{1}{\eta}[x^t - \overline{x}^{t+1}]$$

and use the definition of $x_a$ (from Algorithm 2) in the above bound to obtain the desired result. $\square$

### Proof of Theorem 3

*Proof.* The proof follows from Theorem 6 with $b = 1$ and the parameters used in the theorem statement. $\square$

### B.2 Proof of Theorem 4

*Proof.* The proof follows from Theorem 6 with $b = n^{2/3}$ and the parameters used in the theorem statement. $\square$

## C Convergence Analysis of PL-variants

### C.1 Proof of Theorem 7

*Proof.* The proof follows immediately from Theorem 8 and Theorem 9 with $b = n^{2/3}$ and the parameters used in theorem statement. $\square$

## C.2 PL-SVRG Convergence Analysis

**Theorem 8.** *Suppose $F$ is a $\mu$-PL function. Let $b \leq n$, $\eta = \rho/L$, $T = \lceil 6L/\rho\mu \rceil$, where $\rho \leq 1/5$ and satisfies the condition*

$$\frac{4\rho^2 m^2}{b} + \rho \leq 1.$$

*Then for the output $x^K$ of PL-SVRG, the following holds:*

$$\mathbb{E}[F(x^K) - F(x^*)] \leq \frac{[F(x^0) - F(x^*)]}{2^K},$$

*where $x^*$ is an optimal solution of Problem* (1).

*Proof.* We define the following:

$$\overline{x}_{t+1}^{s+1} = \text{prox}_{\eta h}(x_t^{s+1} - \eta \nabla f(x_t^{s+1})). \tag{31}$$

We first analyze one iteration of the PROXSVRG for PL functions. PL-SVRG essentially uses this as subroutine multiple times in order to obtain the final solution. The proof is similar to that of Theorem 4 until Equation (11). We have the following inequalities:

$$\mathbb{E}[F(\overline{x}_{t+1}^{s+1})] \leq \mathbb{E}\left[F(x_t^{s+1}) + \left[\frac{L}{2} - \frac{1}{\eta}\right]\|\overline{x}_{t+1}^{s+1} - x_t^{s+1}\|^2\right], \tag{32}$$

$$\mathbb{E}[F(x_{t+1}^{s+1})] \leq \mathbb{E}\Bigg[F(\overline{x}_{t+1}^{s+1}) + \langle x_{t+1}^{s+1} - \overline{x}_{t+1}^{s+1}, \nabla f(x_t^{s+1}) - v_t^{s+1}\rangle$$

$$+ \left[\frac{L}{2} - \frac{1}{2\eta}\right]\|x_{t+1}^{s+1} - x_t^{s+1}\|^2 + \left[\frac{L}{2} + \frac{1}{2\eta}\right]\|\overline{x}_{t+1}^{s+1} - x_t^{s+1}\|^2 - \frac{1}{2\eta}\|x_{t+1}^{s+1} - \overline{x}_{t+1}^{s+1}\|^2\Bigg]. \tag{33}$$

Furthermore, we have the following inequality:

$$\mathbb{E}[F(\overline{x}_{t+1}^{s+1})] \leq \mathbb{E}\left[F(x_t^{s+1}) + \langle \nabla f(x_t^{s+1}), \overline{x}_{t+1}^{s+1} - x_t^{s+1}\rangle + \frac{L}{2}\|\overline{x}_{t+1}^{s+1} - x_t^{s+1}\|^2 + h(\overline{x}_{t+1}^{s+1}) - h(x_t^{s+1})\right]$$

$$\leq \mathbb{E}\left[F(x_t^{s+1}) + \langle \nabla f(x_t^{s+1}), \overline{x}_{t+1}^{s+1} - x_t^{s+1}\rangle + \frac{1}{2\eta}\|\overline{x}_{t+1}^{s+1} - x_t^{s+1}\|^2 + h(\overline{x}_{t+1}^{s+1}) - h(x_t^{s+1})\right]$$

$$= \mathbb{E}\left[F(x_t^{s+1}) - \frac{\eta}{2}D_h(x_t^{s+1}, \frac{1}{\eta})\right] \leq \mathbb{E}\left[F(x_t^{s+1}) - \frac{\eta}{2}D_h(x_t^{s+1}, \mu)\right]$$

$$\leq \mathbb{E}\left[F(x_t^{s+1}) - \mu\eta[F(x_t^{s+1}) - F(x^*)]\right] \tag{34}$$

The first inequality follows from Lipschitz continuity of the gradient of $f$. The second inequality follows from the fact that $\eta < 1/L$. The third inequality follows from the fact that $D_h(x,.)$ is a decreasing function. Here, we are implicitly using the fact that $\mu \leq L$ (which can be shown easily for $\mu$-PL functions that are $L$-smooth). Adding 2/3× Equation (32) and 1/3× Equation (34), we have the following:

$$\mathbb{E}[F(\overline{x}_{t+1}^{s+1})] \leq \mathbb{E}\left[F(x_t^{s+1}) + \left[\frac{L}{3} - \frac{2}{3\eta}\right]\|\overline{x}_{t+1}^{s+1} - x_t^{s+1}\|^2 - \frac{\mu\eta}{3}[F(x_t^{s+1}) - F(x^*)]\right]. \tag{35}$$

Adding the above equation with Equation (33), we have the following:

$$\mathbb{E}[F(x_{t+1}^{s+1})] \leq \mathbb{E}\Bigg[F(x_t^{s+1}) + \left[\frac{5L}{6} - \frac{1}{6\eta}\right]\|\overline{x}_{t+1}^{s+1} - x_t^{s+1}\|^2 + \left[\frac{L}{2} - \frac{1}{2\eta}\right]\|x_{t+1}^{s+1} - x_t^{s+1}\|^2$$

$$- \frac{\mu\eta}{3}[F(x_t^{s+1}) - F(x^*)] - \frac{1}{2\eta}\|x_{t+1}^{s+1} - \overline{x}_{t+1}^{s+1}\|^2 + \langle x_{t+1}^{s+1} - \overline{x}_{t+1}^{s+1}, \nabla f(x_t^{s+1}) - v_t^{s+1}\rangle\Bigg]. \tag{36}$$

Using Cauchy-Schwarz and Young's inequality and the fact that $\eta \leq 1/5L$, we have the following:

$$\mathbb{E}[F(x_{t+1}^{s+1})] \tag{37}$$

$$\leq \mathbb{E}\left[F(x_t^{s+1}) + \left[\tfrac{L}{2} - \tfrac{1}{2\eta}\right]\|x_{t+1}^{s+1} - x_t^{s+1}\|^2 - \tfrac{\mu\eta}{3}[F(x_t^{s+1}) - F(x^*)] + \tfrac{\eta}{2}\|\nabla f(x_t^{s+1}) - v_t^{s+1}\|^2\right]$$

$$\leq \mathbb{E}\left[F(x_t^{s+1}) + \left[\tfrac{L}{2} - \tfrac{1}{2\eta}\right]\|x_{t+1}^{s+1} - x_t^{s+1}\|^2 - \tfrac{\mu\eta}{3}[F(x_t^{s+1}) - F(x^*)] + \tfrac{\eta L^2}{2b}\|x_t^{s+1} - \tilde{x}^s\|^2\right]. \tag{38}$$

The second inequality follows from Lemma 3. We use the similar proof technique as in Theorem 5 and 6 and define the following lyapunov function: $R_{t+1}^{s+1} = \mathbb{E}[F(x_{t+1}^{s+1}) + c_{t+1}\|x_{t+1}^{s+1} - \tilde{x}^s\|^2]$. Let $\beta = b/n$. Using the bound on the lyapunov function in Equation (15), we have the following:

$$R_{t+1}^{s+1} \leq \mathbb{E}\Big[F(x_t^{s+1}) - \tfrac{\mu\eta}{3}[F(x_t^{s+1}) - F(x^*)] + \left[c_{t+1}\left(1 + \tfrac{1}{\beta}\right) + \tfrac{L}{2} - \tfrac{1}{2\eta}\right]\|x_{t+1}^{s+1} - x_t^{s+1}\|^2$$
$$\qquad + \left[c_{t+1}(1+\beta) + \tfrac{\eta L^2}{2b}\right]\|x_t^{s+1} - \tilde{x}^s\|^2\Big]$$
$$\leq \mathbb{E}\left[F(x_t^{s+1}) - \tfrac{\mu\eta}{3}[F(x_t^{s+1}) - F(x^*)] + \left[c_{t+1}(1+\beta) + \tfrac{\eta L^2}{2b}\right]\|x_t^{s+1} - \tilde{x}^s\|^2\right]$$
$$= R_t^{s+1} - \tfrac{\mu\eta}{3}\mathbb{E}[F(x_t^{s+1}) - F(x^*)]. \tag{39}$$

The second inequality follows from the fact that:

$$c_{t+1}\left(1 + \frac{1}{\beta}\right) + \frac{L}{2} \leq \frac{1}{2\eta}.$$

This, again, follows from argument stated in Theorem 5, the fact that $\eta = \rho/L$ and

$$\frac{4\rho^2 m^2}{b} + \rho \leq 1.$$

Adding Equation (39) across all the iterations epoch $s+1$ and then using telescopy sum, we get

$$R_m^{s+1} \leq R_0^{s+1} - \sum_{t=0}^{m-1}\frac{\mu\eta}{3}\mathbb{E}[F(x_t^{s+1}) - F(x^*)]. \tag{40}$$

We observe that $R_m^{s+1} = \mathbb{E}[F(x_m^{s+1})] = \mathbb{E}[F(\tilde{x}^{s+1})]$. This is due the fact that $c_m = 0$ and the definition of $\tilde{x}^{s+1}$. Furthermore, $R_0^{s+1} = \mathbb{E}[F(x_0^{s+1})] = \mathbb{E}[F(\tilde{x}^s)]$. This is due to the fact that $x_0^{s+1} = \tilde{x}^s$. Therefore, using the above reasoning in Equation (40), we have

$$\mathbb{E}[F(\tilde{x}^{s+1})] \leq \mathbb{E}[F(\tilde{x}^s)] - \sum_{t=0}^{m-1}\frac{\mu\eta}{3}\mathbb{E}[F(x_t^{s+1}) - F(x^*)].$$

Adding the inequality stated above across all the epochs and using telescopy sum, we have:

$$\sum_{s=0}^{S}\sum_{t=0}^{m-1}\frac{\mu\eta}{3}\mathbb{E}[F(x_t^{s+1}) - F(x^*)] \leq \mathbb{E}[F(x^0)] - \mathbb{E}[F(\tilde{x}^S)] \leq \mathbb{E}[F(x^0)] - F(x^*).$$

The second inequality follows from the optimality of $x^*$. Using the definition of $x^k$ in PL-SVRG, we have the following:

$$\mathbb{E}[F(x^1) - F(x^*)] \leq \frac{3\mathbb{E}[F(x^0) - F(x^*)]}{\mu\eta T}$$
$$\leq \frac{\mathbb{E}[F(x^0) - F(x^*)]}{2}.$$

The second inequality follows from the fact that $T = \lceil 6L/\rho\mu \rceil$. Using this recursion, we have the desired result. $\square$

## C.3 PL-SAGA Convergence Analysis

**Theorem 9.** *Suppose $F$ is a $\mu$-PL function. Let $b \leq n$, $\eta = \rho/L$, $T = \lceil 6L/\mu\rho \rceil$ where $\rho \leq 1/5$ and satisfies the following condition:*
$$\frac{16n^2\rho^2}{b^3} + \rho \leq 1.$$

*Then for the output $x^K$ of* PL-SAGA*, the following holds:*
$$\mathbb{E}[F(x^K) - F(x^*)] \leq \frac{[F(x^0) - F(x^*)]}{2^K},$$

*where $x^*$ is an optimal solution of Problem* (1).

*Proof.* We define the following:
$$\overline{x}^{t+1} = \text{prox}_{\eta h}(x^t - \eta \nabla f(x^t)). \tag{41}$$

Similar to Theorem 8, we start with one iteration of PL-SAGA algorithm. In particular, we first analyze the case of $T$ iterations of SAGA. Further recursing on the the result obtain will give us the desired result. The first part of the theorem is similar to the proof in Theorem 8. Using essentially a similar argument as the one in Theorem 8 until Equation (32), we have the following:

$$\mathbb{E}[F(x^{t+1})] \leq \mathbb{E}\left[F(x^t) + \left[\frac{L}{2} - \frac{1}{2\eta}\right] \|x^{t+1} - x^t\|^2 - \frac{\mu\eta}{3}[F(x^t) - F(x^*)] + \frac{\eta L^2}{2nb} \sum_{i=1}^{n} \|x^t - \alpha_i^t\|^2\right]. \tag{42}$$

We the following Lyapunov function:
$$R_t = \mathbb{E}[F(x^t) + \frac{c_t}{n} \sum_{i=1}^{n} \|x^t - \alpha_i^t\|^2],$$

as defined in Theorem 6. Using the same argument in Theorem 6 to bound it, we have the following:

$$R_{t+1} \leq \mathbb{E}\Big[F(x^t) - \frac{\mu\eta}{3}[F(x^t) - F(x^*)] + \left[c_{t+1}\left(1 + \frac{1-p}{\beta}\right) + \frac{L}{2} - \frac{1}{2\eta}\right] \|x^{t+1} - x^t\|^2$$
$$+ \left[\frac{c_{t+1}(1+\beta)(1-p)}{n} + \frac{\eta L^2}{2nb}\right] \sum_{i=1}^{n} \|x^t - \alpha_i^t\|^2\Big]$$
$$\leq R_t - \frac{\mu\eta}{3}\mathbb{E}[F(x^t) - F(x^*)]. \tag{43}$$

Recall that $p = 1 - (1 - 1/n)^b$. The second inequality is due to the following inequality:
$$c_{t+1}\left(1 + \frac{1-p}{\beta}\right) + \frac{L}{2} \leq \frac{1}{2\eta}.$$

This is obtained by the same argument in Theorem 6. Adding Equation (43) over all the iterations and using telescopy sum, we have the following:
$$\mathbb{E}[F(x^T)] \leq \mathbb{E}[F(x^0)] - \sum_{t=0}^{T-1} \frac{\mu\eta}{3}\mathbb{E}[F(x^t) - F(x^*)].$$

The above inequality is obtained from the fact that $R_T = \mathbb{E}[F(x^T)]$. This is due the fact that $c_T = 0$. Furthermore, $R_0 = \mathbb{E}[F(x^0)]$. This is due to the fact that $\alpha_i^0 = x^0$ for all $i \in [n]$. Therefore, we have:
$$\sum_{t=0}^{T-1} \frac{\mu\eta}{3}\mathbb{E}[F(x^t) - F(x^*)] \leq \mathbb{E}[F(x^0)] - \mathbb{E}[F(x^T)] \leq \mathbb{E}[F(x^0) - F(x^*)].$$

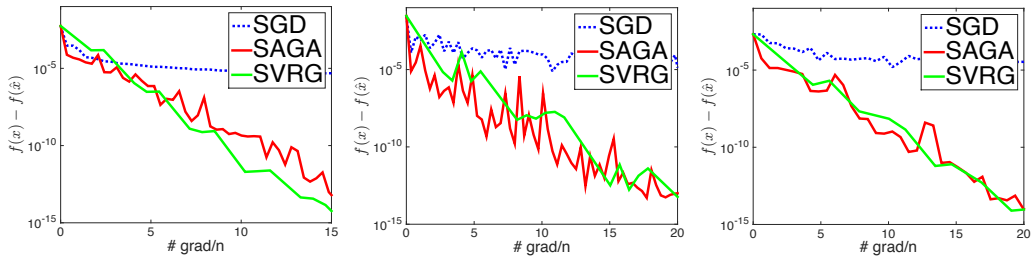

Figure 3: Non-negative principal component analysis. Performance of SGD, PROXSVRG and PROXSAGA on 'real-sim' (left), 'covtype'(center) and 'ijcnn1' (right) datasets. Recall that the y-axis is the function suboptimality i.e., $f(x) - f(\hat{x})$ where $\hat{x}$ represents the best solution obtained by running gradient descent for long time and with multiple restarts.

Using the definition of $x^k$ in PL-SAGA, we have the following:

$$\mathbb{E}[F(x^1) - F(x^*)] \leq \frac{3\mathbb{E}[F(x^0) - F(x^*)]}{\mu\eta T}$$
$$\leq \frac{\mathbb{E}[F(x^0) - F(x^*)]}{2}.$$

The second inequality follows from the fact that $T = \lceil 6L/\mu\rho \rceil$. Using the above recursion repeatedly, we obtain the desired result. □

## D   Additional Experiments

We present the additional experiments for non-negative principal component analysis problems in this section. Figure 3 shows the additional results. Similar to Figure 2, we see that PROXSVRG and PROXSAGA outperform PROXSGD. We did not find any significant difference in the performance of PROXSVRG and PROXSAGA.

## E   Lemmatta

We first few intermediate results that are useful for our analysis. These results are key to the mirror descent analysis [15]. We prove them here for completeness.

**Lemma 1.** *Suppose we define the following:*

$$y = prox_{\eta h}(x - \eta d'). \tag{44}$$

*for some $d' \in \mathbb{R}^d$. Then for $y$, the following inequality holds:*

$$h(y) + \langle y - z, d' \rangle \leq h(z) + \tfrac{1}{2\eta}\left[\|z - x\|^2 - \|y - x\|^2 - \|y - z\|^2\right]. \tag{45}$$

*for all $z \in \mathbb{R}^d$.*

*Proof.* From Lemma 5 applied on Equation (44), we get the following:

$$h(y) + \langle y - x, d' \rangle + \tfrac{1}{2\eta}\|y - x\|^2 + \frac{\eta}{2}\|d'\|^2$$
$$= h(y) + \tfrac{1}{2\eta}\|y - (x - \eta d')\|^2$$
$$\leq h(z) + \tfrac{1}{2\eta}\|z - (x - \eta d')\|^2 - \tfrac{1}{2\eta}\|y - z\|^2$$
$$= h(z) + \langle z - x, d' \rangle + \tfrac{1}{2\eta}\|z - x\|^2 + \frac{\eta}{2}\|d'\|^2 - \tfrac{1}{2\eta}\|y - z\|^2. \tag{46}$$

By rearranging Equation (46), we obtain the following inequality that concludes the proof.

$$h(y) + \langle y - z, d' \rangle \leq h(z) + \tfrac{1}{2\eta}\left[\|z - x\|^2 - \|y - x\|^2 - \|y - z\|^2\right].$$

□

The following key lemma involving function $F$ is useful for proving the convergence of PROXSVRG and PROXSAGA.

**Lemma 2.** *Suppose we define the following:*

$$y = \text{prox}_{\eta h}(x - \eta d').$$

*for some $d' \in \mathbb{R}^d$. Then for $y$, the following inequality holds:*

$$F(y) \leq F(z) + \langle y - z, \nabla f(x) - d' \rangle$$
$$+ \left[\tfrac{L}{2} - \tfrac{1}{2\eta}\right]\|y - x\|^2 + \left[\tfrac{L}{2} + \tfrac{1}{2\eta}\right]\|z - x\|^2 - \tfrac{1}{2\eta}\|y - z\|^2. \tag{47}$$

*for all $z \in \mathbb{R}^d$.*

*Proof.* We have the following inequalities for function $f$:

$$f(y) \leq f(x) + \langle \nabla f(x), y - x \rangle + \frac{L}{2}\|y - x\|^2,$$
$$f(x) \leq f(z) + \langle \nabla f(x), x - z \rangle + \frac{L}{2}\|x - z\|^2.$$

The above inequalities can be obtained by application of Lemma 6. Adding both the inequalities above, we obtain the following inequality:

$$f(y) \leq f(z) + \langle \nabla f(x), y - z \rangle + \tfrac{L}{2}\left[\|y - x\|^2 + \|z - x\|^2\right]. \tag{48}$$

Adding Equations (45) (which follows from Lemma 1) and (48), we obtain the inequality:

$$F(y) \leq F(z) + \langle y - z, \nabla f(x) - d' \rangle$$
$$+ \left[\tfrac{L}{2} - \tfrac{1}{2\eta}\right]\|y - x\|^2 + \left[\tfrac{L}{2} + \tfrac{1}{2\eta}\right]\|z - x\|^2 - \tfrac{1}{2\eta}\|y - z\|^2. \tag{49}$$

Here we used that definition $F(x) = f(x) + h(x)$. This concludes our proof. □

The following result is useful for bounding the variance of the updates of PROXSVRG and follows from slight modification of result in [22]. We give the proof here for completeness.

**Lemma 3** ([22]). *For the iterates $x_t^{s+1}, v_t^{s+1}$ and $\tilde{x}^s$ where $t \in \{0, \ldots, m-1\}$ and $s \in \{0, \ldots, S-1\}$ in Algorithm 1, the following inequality holds:*

$$\mathbb{E}[\|\nabla f(x_t^{s+1}) - v_t^{s+1}\|^2] \leq \frac{L^2}{b}\|x_t^{s+1} - \tilde{x}^s\|^2.$$

*Proof.* Let us define the following notation for the ease of exposition:

$$\zeta_t^{s+1} = \frac{1}{|I_t|} \sum_{i \in I_t} \left( \nabla f_i(x_t^{s+1}) - \nabla f_i(\tilde{x}^s) \right).$$

Using this notation, we obtain the following bound:

$$\mathbb{E}[\|\nabla f(x_t^{s+1}) - v_t^{s+1}\|^2] = \mathbb{E}[\|\zeta_t^{s+1} + \nabla f(\tilde{x}^s) - \nabla f(x_t^{s+1})\|^2]$$

$$= \mathbb{E}[\|\zeta_t^{s+1} - \mathbb{E}[\zeta_t^{s+1}]\|^2] = \frac{1}{b^2} \mathbb{E}\left[\left\|\sum_{i \in I_t} \left(\nabla f_i(x_t^{s+1}) - \nabla f_i(\tilde{x}^s) - \mathbb{E}[\zeta_t^{s+1}]\right)\right\|^2\right]$$

The second equality is due to the fact that $\mathbb{E}[\zeta_t^{s+1}] = \nabla f(x_t^{s+1}) - \nabla f(\tilde{x}^s)$. From the above relationship, we get

$$\mathbb{E}[\|\nabla f(x_t^{s+1}) - v_t^{s+1}\|^2] \leq \frac{1}{b} \mathbb{E}\left[\sum_{i \in I_t} \|\nabla f_i(x_t^{s+1}) - \nabla f_i(\tilde{x}^s) - \mathbb{E}[\zeta_t^{s+1}]\|^2\right]$$

$$\leq \frac{1}{b} \mathbb{E}\left[\sum_{i \in I_t} \|\nabla f_i(x_t^{s+1}) - \nabla f_i(\tilde{x}^s)\|^2\right] \leq \frac{L^2}{b} \mathbb{E}[\|x_t^{s+1} - \tilde{x}^s\|^2]$$

The first inequality follows from Lemma 7. The second inequality is due to the fact that for a random variable $\zeta$, $\mathbb{E}[\|\zeta - \mathbb{E}[\zeta]\|^2] \leq \mathbb{E}[\|\zeta\|^2]$. The last inequality follows from $L$-smoothness of $f_i$. □

A similar result can be obtained for PROXSAGA. The key difference with that of Lemma 3 is that the variance term in PROXSAGA involves $\alpha_i^t$. Again, we provide the proof for completeness.

**Lemma 4.** *For the iterates $x^t, v^t$ and $\{\alpha_i^t\}_{i=1}^n$ where $t \in \{0, \ldots, T-1\}$ in Algorithm 2, the following inequality holds:*

$$\mathbb{E}[\|\nabla f(x^t) - v^t\|^2] \leq \frac{L^2}{nb} \sum_{i=1}^n \mathbb{E}\|x^t - \alpha_i^t\|^2.$$

*Proof.* Let us define the following notation for the ease of exposition:

$$\zeta_t = \frac{1}{|I_t|} \sum_{i \in I_t} \left( \nabla f_i(x^t) - \nabla f_i(\alpha_i^t) \right).$$

In this notation, we have the following:

$$\mathbb{E}[\|\nabla f(x^t) - v^t\|^2] = \mathbb{E}\left[\left\|\zeta_t + \frac{1}{n} \sum_{i=1}^n \nabla f(\alpha_i^t) - \nabla f(x^t)\right\|^2\right]$$

$$= \mathbb{E}[\|\zeta_t - \mathbb{E}[\zeta_t]\|^2] = \frac{1}{b^2} \mathbb{E}\left[\left\|\sum_{i \in I_t} \left(\nabla f_i(x^t) - \nabla f_i(\alpha_i^t) - \mathbb{E}[\zeta_t]\right)\right\|^2\right]$$

The second equality follows from the fact that $\mathbb{E}[\zeta_t] = \nabla f(x^t) - \frac{1}{n}\sum_{i=1}^n \nabla f(\alpha_i^t)$. From the above inequality, we get

$$\mathbb{E}[\|\nabla f(x^t) - v^t\|^2] \leq \frac{1}{b} \mathbb{E}\left[\sum_{i \in I_t} \|\nabla f_i(x^t) - \nabla f_i(\alpha_i^t) - \mathbb{E}[\zeta_t]\|^2\right]$$

$$\leq \frac{1}{b} \mathbb{E}\left[\sum_{i \in I_t} \|\nabla f_i(x^t) - \nabla f_i(\alpha_i^t)\|^2\right] \leq \frac{L^2}{nb} \sum_{i=1}^n \mathbb{E}[\|x^t - \alpha_i^t\|^2]$$

The first inequality is due to Lemma 7. The second inequality follows from the fact that for a random variable $\zeta$, $\mathbb{E}[\|\zeta - \mathbb{E}[\zeta]\|^2] \leq \mathbb{E}[\|\zeta\|^2]$. The last inequality is from $L$-smoothness of $f_i$ for all $i \in [n]$ and uniform randomness of the set $I_t$. □

The following lemma is a classical result in mirror descent analysis.

**Lemma 5.** *Suppose function $h : \mathbb{R}^d \to \mathbb{R}$ is l.s.c and $y = prox_{\eta h}(x)$. Then we have the following:*

$$h(y) + \tfrac{1}{2\eta}\|y - x\|^2 \leq h(z) + \tfrac{1}{2\eta}\|z - x\|^2 - \tfrac{1}{2\eta}\|y - z\|^2,$$

*for all $z \in \mathbb{R}^d$.*

**Lemma 6.** *Suppose function $f : \mathbb{R}^d \to \mathbb{R}$ is $L$-smooth, then we have the following:*

$$f(y) + \langle \nabla f(y), x - y \rangle - \frac{L}{2}\|x - y\|^2 \leq f(x) \leq f(y) + \langle \nabla f(y), x - y \rangle + \frac{L}{2}\|x - y\|^2,$$

*for all $x, y \in \mathbb{R}^d$.*

**Lemma 7.** *For random variables $z_1, \ldots, z_r$ are independent and mean 0, we have*

$$\mathbb{E}\left[\|z_1 + \ldots + z_r\|^2\right] = \mathbb{E}\left[\|z_1\|^2 + \ldots + \|z_r\|^2\right].$$

*Proof.* We have the following:

$$\mathbb{E}\left[\|z_1 + \ldots + z_r\|^2\right] = \sum_{i,i=1}^{r} \mathbb{E}\left[z_i z_j\right] = \mathbb{E}\left[\|z_1\|^2 + \ldots + \|z_r\|^2\right].$$

The second equality follows from the fact that $z_i$'s are independent and mean 0. □

**Lemma 8.** *For random variables $z_1, \ldots, z_r$, we have*

$$\mathbb{E}\left[\|z_1 + \ldots + z_r\|^2\right] \leq r\mathbb{E}\left[\|z_1\|^2 + \ldots + \|z_r\|^2\right].$$